\documentclass[a4paper,12pt]{amsart}
\usepackage{amsmath,amssymb}
\def\C{\Bbb C}
\def\D{\Delta}

\def\Re{\mbox{Re}}
\def\O{\mathcal O}
\def\d{\delta}
\def\e{\eta}
\def\ww{\varphi(\eta)}
\def\z{\zeta}
\def\zz{\varphi(\zeta)}
\def\eps{\varepsilon}

\def\k{\kappa}
\def\vp{\varphi}
\def\bd{\partial}
\def\td{\tilde\delta}

\def\FF{\mathcal F_c}
\def\ov{\overline}

\renewcommand\d{\delta}

\newcommand\CC{{\mathbb C}}

\newcommand\RR{{\mathbb R}}
\newcommand\DD{{\mathbb D}}
\newcommand\TT{{\mathbb T}}

\def\re{\operatorname{Re}}
\def\diam{\operatorname{diam}}

\newtheorem{thm}{Theorem}
\newtheorem{prop}[thm]{Proposition}
\newtheorem{cor}[thm]{Corollary}
\theoremstyle{remark}\newtheorem{remark}{Remark}

\title[The Kobayashi distance and complex geodesics]
{Lower estimates of the Kobayashi distance and limits of complex geodesics}

\author[\L ukasz Kosi\'nski]{\L ukasz Kosi\'nski}
\address{\L. Kosi\'nski: Institute of Mathematics, Faculty of Mathematics and Computer Science,
	Jagiellonian University, \L ojasiewicza 6, 30-348 Krak\'ow, Poland} \email{lukasz.kosinski@uj.edu.pl}

\author{Nikolai Nikolov}
\address{N. Nikolov\\Institute of Mathematics and Informatics\\Bulgarian Academy
	of Sciences\\Acad. G. Bonchev 8, 1113 Sofia, Bulgaria
	\vspace{1mm}
	\newline Faculty of Information Sciences\\
	State University of Library Studies and Information Technologies\\
	Shipchenski prohod 69A, 1574 Sofia,
	Bulgaria}

\email{nik@math.bas.bg}

\thanks{The first named author was supported by the NCN grant SONATA BIS no. 2017/26/E/ST1/00723.
	The second named author was partially supported by the Bulgarian National Science Fund,
	Ministry of Education and Science of Bulgaria under contract KP-06-N52/3.}


\begin{document}
	
	
	\begin{abstract} {It is proved for a strongly pseudoconvex domain $D$ in $\C^d$ with $\mathcal C^{2,\alpha}$-smooth boundary that any complex geodesic through every two close points of $D$ sufficiently close to $\bd D$ and whose difference is non-tangential to $\bd D$ intersect a compact subset of $D$ that depends only on the rate of non-tangentiality.~As an application, a lower bound for the Kobayashi distance is obtained.}
	\end{abstract}
	
	\maketitle
	
	\section{Introduction}
	
In \cite{H2} Huang proved the localization theorem for infinitesimal extremal mappings near $\mathcal C^3$-smooth strongly pseudoconvex points.

This result, which was initially developed to solve a conjecture of
Abate-Vensentini in \cite{H1}, has found many other applications in working on various problems. For instance, it plays a crucial role in solving the homogeneous Monge-Amp\`ere equations
or in establishing various regularity for CR mappings.
Its immediate consequence was also the existence of geodesics with prescribed boundary data, a result previously obtained in \cite{CHL} within the class of $\mathcal C^{14}$-smooth strictly linearly convex domains.
Recently, Huang \cite{H3} showed that his result remains true for strongly pseudoconvex point with only $\mathcal C^{2, \alpha}$-regularity, $\alpha\in (0,1]$.

In this paper we extend the localization theorem to non-infinitesimal settings, while the estimate we are able to achieve applied to infinitesimal case is sharper.

Some part of our paper is devoted to yet another application: we obtain lower bounds for the Kobayashi distance. Estimates that we are going to present being in the spirit of Balogh-Bonk results \cite{BB} are sensitive in the non-tangential case.

	\section{Definitions and main results}

Let $D$ be a domain in $\C^d$, $z,w\in D,$ and $X\in\C^d.$ The Kobayashi distance $k_D$ is the largest pseudodistance not exceeding the Lempert function
$$l_D(z,w)=\inf\{\tanh^{-1}|\alpha|:\exists\vp\in\O(\D,D)
	\hbox{ with }\varphi(0)=z,\varphi(\alpha)=w\},$$
where $\D$ is the unit disc and $\tanh^{-1}t=\frac12\log\frac{1+t}{1-t}.$	This distance is the integrated form of the Kobayashi metric
$$\k_D(z;X)=\inf\{|\alpha|:\exists\vp\in\mathcal O(\D,D), \vp(0)=z, \alpha\vp'(0)=X\}.$$
	
	Set $\d_D(z)=\mbox{dist}(z,\bd D)$ and $h_D(z,w)=\d_D(z)^{1/2}\d_D(w)^{1/2}.$
	
	Recall that if $p$ is a Dini-smooth boundary point of $D$, then
	(see \cite[Corollary 8]{NA})
	\begin{equation}\label{na}
	k_D(z,w)\le\log\left(1+\frac{2|z-w|}{h_D(z,w)}\right),\quad z,w\mbox{ near }p.
	\end{equation}
	
	The example of the unit ball shows a similar lower bound has to include the term
	$\frac{|z-w|^2}{\d_D(z)^{1/2}\d_D(w)^{1/2}}.$
	
	Let us mention two results in this direction. If $p$ is a strongly pseudoconvex boundary point of $D$, then
	there is $c>0$ such that (see \cite[Theorems 1.6 and 1.7]{NT})
	\begin{equation}\label{nt}
	k_D (z,w)\ge\log\left(1+\frac{c|z-w|}{\d_D(z)^{1/2}}\right)
	\left(1+\frac{c|z-w|}{\d_D(w)^{1/2}}\right),\quad z,w\mbox{ near }p.
	\end{equation}
	
For strongly pseudoconvex domain $D,$ Balogh and Bonk introduced a positive function $g_D$ on $D\times D$ which depends on $\d_D$ and the Carnot-Carath\'eodory metric on $\partial D$ (see  \cite[formula (1.2)]{BB}) and they proved that $g_D-k_D$ is bounded on $D\times D.$ As a consequence of their work, there exists $C>0$ such that
	\begin{equation}\label{bb}
	\log\left(1+\frac{|z-w|^2}{\d_D(z)^{1/2}\d_D(w)^{1/2}}\right)-C\le k_D (z,w)
	\end{equation}
$$\le\log\left(1+\frac{|z-w|}{\d_D(z)^{1/2}\d_D(w)^{1/2}}\right)+C,$$
	which also follows from \eqref{na} and \eqref{nt}.
	
The following estimate is more sensitive in the non-tangential case.
	
\begin{prop}\label{c2} Let $p$ be a strongly pseudoconvex boundary point of a domain $D.$
Then there exist a neighborhood $U$ of $p$ and a constant $c>0$ such that
$$k_D(z,w)\ge\log\left(1+\frac{c|(z-w)_n|}{\d_D(z)^{1/2}\d_D(w)^{1/2}}\right),\quad z,w\in D\cap U.$$
	\end{prop}
	
	Here $X_n$ denotes the projection of $X$ on the real normal to $\bd D$ at the closest point $p(z)\in\bd D$
	to $z.$ Note that $|X_n|=2|\Re\langle X,g_D(z)\rangle|,$
	where $g_D(z)=:\bar\partial\td_D(z)=\bar\partial\td_D(p(z))$ and $\td_D$ is the signed distance to $\partial D.$
	\smallskip
	
Proposition \ref{c2} together with \eqref{nt} leads to the following more precise estimate:
	$$k_D(z,w)\ge\log\left(1+c\left(\frac{|(z-w)_n|+|z-w|^2}{\d_D(z)^{1/2}\d_D(w)^{1/2}}+\frac{|z-w|}{\d_D(z)^{1/2}}+
	\frac{|z-w|}{\d_D(w)^{1/2}}\right)\right).$$ Actually, it is a consequence of the trivial inequality: $\max\{ \log(1+ x_1), \log(1 + x_2)\} \geq \log\left( 1 + \frac{x_1 + x_2}{2}\right).$
	
	One may expect similar results, replacing the real by the complex normal.
	
	\begin{prop}\label{c3} Let $p$ be a $\mathcal C^{2,\alpha}$-smooth strongly pseudoconvex boundary point of a domain $D.$
		Then there exists a neighborhood $U$ of $p$ such that for any $\eps\in(0,1]$ one may find a constant $c=c(\eps)>0$ such that
		if $z,w\in D\cap U$ and $|(z-w)_N|\ge\eps|z-w|,$ then
		$$k_D(z,w)\ge\log\left(1+\frac{c|(z-w)_N|}{\d_D(z)^{1/2}\d_D(w)^{1/2}}\right).$$
	\end{prop}
	
	Here and in the sequel $\alpha\in (0,1]$ and $X_N$ denotes the projection of $X$ on the complex normal to $\bd D$
    at the closest point $p(z)\in\bd D$	to $z.$ Note that $|X_N|=2|\langle X,g_D(z)\rangle|.$
	
	It is natural to ask if this result remains true for $\eps=0,$ that is, without the assumption about $(z-w)_N.$
	Proposition \ref{c33} below gives a partial answer to this question.
	
Thy key point in the proof of Proposition \ref{c3} is the following result which is also of independent interest. To formulate it let us introduce a class of {\it $FH$-domains}. It is composed of bounded pseudoconvex domains that either have $\mathcal C^\infty$-smooth boundary, or their closures posses Stein neighborhood bases. The motivation for the name that we are using here comes from the exposing type theorem that is essentially due to the work of Huang (see \cite[Lemma 1]{H1}) for applying a local version  of the Fornaess Forn{\ae}ss embedding theorem \cite{For}. Note that any strongly pseudoconvex domain is an $FH$-domain.

	\begin{thm}\label{th:gen} Let $D$ be an $FH$-domain. Suppose that $p$ is a $\mathcal C^{2, \alpha}$-smooth strongly pseudoconvex boundary point of $D$. Then there is a neighborhood $U$ of $p$ and a constant $C>0$ such
		$$|(z-w)_N| \leq C |z-w| \diam (\varphi)$$ for any extremal mapping $\varphi$ through $z$ and $w$, where $z\neq w\in D\cap U$.
	\end{thm}

As we mentioned in the introduction the infinitesimal counterpart of the above result has been obtained in \cite{H2, H3}. Note that it can be deduced directly from Theorem~\ref{th:gen}, as the proof of \cite[Theorem 1]{H1} implies that any extremal mapping whose diameter is small enough must be a complex geodesic.

The proof presented here is based on the deep theory of Lempert \cite{L1,L2} and its extension to $\mathcal C^{2,\alpha}$-case due to Huang and Lempert. Instead of considering additional Lagrange-Euler equations or investigating subtle (e.g. boundary) perturbation of stationary discs, the main point of our argument relies on a careful analysis of the scaling carried out in \cite{L1}.
	
	Theorem~\ref{th:gen} together with \cite[Proposition 4]{CHL} (see also \eqref{d2} below) leads to the following
	\begin{cor}\label{geod} Let $D$ and $p$ be as in Theorem~\ref{th:gen}.
		Then for any $\eps\in(0,1]$ there exist a neighborhood $U$ of $p$ and a compact subset $K$ of $D$ such that
		if $z\neq w\in D\cap U$ and $|(z-w)_N|\ge\eps|z-w|,$ then any extremal mapping through $z$ and $w$ intersects $K$.
	\end{cor}
	
	The respective definitions, as well as other facts about complex geodesic, are given in section \ref{cg}.

	\section{Extremal mappings and complex geodesics}\label{cg}
	
	Let $D$ be a domain in $\C^d,$ $z\neq w\in D,$ and $0\neq X\in\C^d.$
	
	A mapping $\vp\in\O(\D,D)$ is called:
	
	-- extremal for $\k_D(z;X):=\alpha$ (in short, $\vp=\vp_{z;X}$)
	if $\vp(0)=z$ and  $\alpha\vp'(0)=X.$

	-- extremal for $k_D(z;w):=\alpha$ (in short, $\vp=\vp_{z,w}$)
	if $\vp(0)=z,$ $\vp(\alpha)=w$
	and $k_D(z,w)=\tanh^{-1}\alpha.$
	
	-- complex geodesic if $k_D(\zz,\ww)=k_\D(\z,\e)$ for any $\z,\e\in\D.$
	\smallskip
	
	Recall that a bounded $D$ is said to be strictly linearly convex if $D$ has $\mathcal C^2$-smooth boundary and, for any $p\in\partial D,$ the restriction of the Hessian of $\td_D$	on $T_p^{\CC}(\bd D)$ is a positive Hermitian form.
	
According to Lempert's theory, if $D$ is bounded strictly linearly convex, there exist unique extremal mappings for $\k_D(z;X)$ and $k_D(z,w),$ and they are complex geodesics. Conversely, after a reparametrization,
	any complex geodesic $\vp$ is extremal for $\k_D(\vp(\z);\vp'(\z))$ and
	$k_D(\vp(\z),\vp(\e))$ for any $\z,\e\in\D.$ It is also known (see \cite{L3}) that if $D$ has a $\mathcal C^{2,\alpha}$-smooth boundary, then
    complex geodesics
	extend $\mathcal C^{1,\alpha-}$-smoothly to $\ov\D$ (that is $\mathcal C^{1,\alpha-\eps}$-smoothly for any $\eps>0$).
	\smallskip
	
	Let us mention the following tangential version of Theorem \ref{th:gen}.
	
	\begin{thm}\label{h1} (\cite[Theorem 1]{H1})
		Let $D$ be an $FH$-domain.
		Suppose that $p$ is a $\mathcal C^{2,\alpha}$-smooth strongly pseudoconvex boundary point of $D.$
		Then for any neighborhood $U$ of $p,$ there exists $\eps>0$ such that
		for each extremal mapping $\vp_{z;X}$ with $|z-p|<\eps$ and
		$|X_N|<\eps|X|,$ one has that $\vp_{z;X}(\D)\subset U$ and $\vp_{z;X}$ is a complex geodesic.
		
		(\cite[Proposition 2.5]{BFW}, see also \cite{Kos}) The same remains true for each extremal mapping $\vp_{z,w}$
		with $|z-p|<\eps, |w-p|<\eps$ and $|z-w|_N<\eps|z-w|.$
	\end{thm}
	
	Let us remark that Theorem~\ref{h1} is stated in \cite{H1} under an additional assumption of $\mathcal C^3$-smoothness, nevertheless the whole proof works in $\mathcal C^{2,\alpha}$-settings, as well.
	
	The next result can be considered as an infinitesimal version of Theorem \ref{th:gen}.
	
	\begin{thm}\label{h2} (\cite[Corollary 1.2]{H3}; \cite[Theorem 2]{H2} if $\alpha=1$)
		Let $p$ be a $\mathcal C^{2,\alpha}$-smooth strongly pseudoconvex point of a domain $D,$ where
		$\alpha\in(0,1].$ Then there exist a neighborhood $U$ of $p$ and $C>0$ such that for each extremal
		mapping $\vp_{z;X}$ with $\vp_{z;X}(\D)\subset U$ one has that $|X_N|<C\diam^\alpha(\vp_{z;X}(\D))|X|.$
	\end{thm}
	
As we already said one can deduce Theorem~\ref{h2} directly from Theorem~\ref{th:gen} with a bit better estimate $|X_N|< C \diam(\varphi_{z;X}(\Delta))|X|$. Note that this estimate is almost precise.
	Indeed, if $D$ is a ball, then
	$$\diam(\vp_{z;X}(\D))\sim \d_D(z)^{1/2}+|X_N|/|X|.\ \footnote
{As usual, $f\lesssim g$ means that $f\le Cg$ for some constant $C>0$ depending only on $D;$ $f\sim g$ means that
$f\lesssim g\lesssim f.$}$$
	This observation can be partially extended in a more general setting.
	
	\begin{prop}\label{diam} Let $D$ be a strongly pseudoconvex domain with
		$\mathcal C^{2,\alpha}$-smooth boundary.
		There exist $C>c>0$ such that for any complex geodesic $\vp\in\O(\D,D)$ one has that
		\begin{equation}\label{d1}
		c\left(\d_D(\zz)^{1/2}+\frac{|\vp'(\z)|_N}{|\vp'(\z)|}\right)\le\diam(\vp),\quad \z\in\D,
		\end{equation}
		\begin{equation}\label{d2}
		c\max(\d_D\circ\vp)^{1/2}\le\diam(\vp)\le C\max(s\circ\d_D\circ\vp),
		\end{equation}
		where $s(x)=-x^{1/2}\log x.$
	\end{prop}
    In \eqref{d1} we may assume that $\diam(\vp)$ is small enough; therefore $\vp(\D)$
	is near some $p\in\partial D$ and hence $|\vp'(\z)|_N$ is well-defined.
	
	One may expect that
	\begin{equation}\label{conj}
	\diam(\vp)\sim\max(\d_D\circ\vp)^{1/2}\sim\max|\vp'|
	\end{equation}
	for any complex geodesic $\vp\in\O(\D,D)$ which is parameterized such that
	$\d_D(\vp(0))=\max(\d_D\circ\vp).$

	\begin{prop}\label{c33}
	Proposition~\ref{c3} remains true for $\eps=0$ if \eqref{conj} holds.
	\end{prop}
	
	Since (cf. \cite[Theorem 19.4.2]{JP})
	\begin{equation}\label{ma}
	\frac{1}{1-|\z|^2}=\k_\D(\z;e)=\k_D(\zz;\vp'(\z))\sim\frac{|\vp'(\z)|_N}{\d_D(\zz)}+\frac{|\vp'(\z)|}{\d_D(\zz)^{1/2}},
	\end{equation}
then
$$\frac{\d_D(\zz)}{\d_{\Delta}(\z)|\vp'(\z)|}\sim\d_D(\zz)^{1/2}+\frac{|\vp'(\z)|_N}{|\vp'(\z)|}$$
	and the inequality \eqref{d1} can be read as
	\begin{equation}\label{equiv}
    \d_D(\zz)\lesssim\diam(\vp)|\vp'(\z)|\d_\D(\z).
    \end{equation}

    Note that, by \cite[Theorem 20]{NO},
    $$\frac{\d_D(\vp(0))}{\d_D(\zz)}\lesssim\frac{1}{\d_\D(\z)}.$$

    On the other hand, if \eqref{conj} is true, then \eqref{equiv} implies the opposite inequality
    and hence we may replace $\lesssim$ by $\sim$ in \eqref{d1} and \eqref{equiv}
    (also, $\max|\vp'|\sim|\vp'(\z)|$ for any $\z\in\D$). This will give us a quantitative version of the first part of Theorem~\ref{h1}.

    Further, for $X\not\in T^{\C}_p(\bd D),$ there exists a unique $\lambda_X\in\TT$ such that
	$\langle\nu_p,\lambda_X X\rangle>0,$ where $\langle\cdot,\cdot\rangle$ is the standard Hermitian
	product and $\nu_p$ is the outward normal to $\partial D$ at $p.$
	
	Theorems \ref{th:gen}, \ref{h1} and \ref{h2} imply the following.
	
	\begin{cor}\label{seq} Let $D$ and $p$ be as in Theorem~\ref{th:gen}. Let $(\vp_{z_j;X_j})$ and $(\psi_{z_j;w_j})$
		be sequences of extremal mappings such that $z_j, w_j\to p,$ $X_j\to X,$ and $\frac{z_j-w_j}{|z_j-w_j|}\to X.$
		
		(a) If $X\in T^{\C}_p(\bd D),$ then $\vp_j, \psi_j\to p$ uniformly on $\ov\D.$
		
		(b) If $X\not\in T^{\C}_p(\bd D)$ and the diameters of the extremal mappings are small enough, then one may find subsequences of $(\vp_{z_j;X_j})$ and $(\psi_{z_j;w_j})$
		which tend, after reparametrizations, resp.~to $\vp$ and $\psi$ in the $\mathcal C^1$-norm on $\ov\D,$
		where $\vp$ and $\psi$ are complex geodesics with $\vp(1)=\psi(1)=p$ and $\vp'(1)=\psi'(1)=\lambda_XX.$
	\end{cor}
	
		According to the recent result of \cite{HW}, geodesics $\vp$ and $\psi$ appearing in Corollary~\ref{seq} coincide up to
	an automorphism of $\D$ if $\bd D$ is $\mathcal C^{2,\alpha},$ where $\alpha>1/2$ (see \cite[Theorem~1.1]{HW} and the discussion in \cite[Remark~2.4]{HW}).

	\section{Proofs of Propositions \ref{c2}, \ref{c3}, \ref{diam}, \ref{c33} and Corollary \ref{seq}}\label{c23}
	\noindent{\it Proof of Proposition \ref{c2}.} It is known that there exists a ball $V$ around $p$ such that for
	any $z\in D$ near $p$ one may find a composition $\Phi_z$ (injective on $V$) of three transformation, all depending
	continuously of $p(z)$ and hence of $z,$ namely the translation $\pi(z)\to 0,$ a unitary transformation and a change of the
	variables of the form $\e=\e_z=(\z_1+c_z\z_1^2+P_z(\z'),\z'),$ where $c_z\in\Bbb R$ and $P_z$ is a quadratic polynomial,
	such that $G_z=\Phi_z(D\cap V)$ is a convex domain and $\re \e_1<0$ is the inner normal to $\bd G_z$ at $0.$
	Set $d_z=\d_{G_z}\circ\Phi_z.$ It is easy to check that near $p$ (see below),
	\begin{equation}\label{appr}
	|(d_z(z)-d_z(w)|=(1+o(1))|(z-w)_n|+O(|z-w|^2),
	\end{equation}
	$$d_z=(1+o(1))\d_D.$$
	
Let now $z,w$ be near $p.$ Assume that $d_z(z)\le d_z(w).$ Note that $G_z\subset\Pi_z=\{\e:\re \e_1<0\},$ and that $k_{\Pi_z}(a,b) = K_{\Pi_1}(a_1, b_1),$ where $\Pi_1= \{z_1:\ \re z_1 <0\}\subset \CC$. Therefore
	$$k_{D\cap V}(z,w)\ge k_{\Pi_z}(\Phi_z(z),\Phi_z(w))\ge\frac{1}{2}\log\frac{\d_{\Pi_z}(\Phi_z(w))}{\d_{\Pi_z}(\Phi_z(z))}$$
	$$\ge\frac{1}{2}\log\frac{d_z(w)}{d_z(z)}\ge\log\left(1+\frac{d_z(w)-d_z(z)}{2d_z(z)^{1/2}d_z(w)^{1/2}}\right).$$
	Hence, switching $z$ and $w$, if necessary, we get
	\begin{equation}\label{imd}
	k_{D\cap V}(z,w)\ge\log\left(1+\frac{|d_z(w)-d_z(z)|}{2d_z(z)^{1/2}d_z(w)^{1/2}}\right).
	\end{equation}
	The same remains true in the case $d_z(z)\ge d_z(w)$ by considering the supporting hyperplane to $\bd G_z$
	at the closest point to $\Phi_z(w).$
	\smallskip
	
On the other hand, the proofs of \cite[Theorem 1.4]{BNT} and \cite[Theorem 1.6]{NT} imply that for any
	neighborhood $U\Subset V$ there exists a constant $C>0$ such that
	\begin{equation}\label{loc}
	k_{D\cap V}(z,w)\le k_D(z,w)+C,\quad k_{D\cap V}\le Ck_D(z,w),\quad z,w\in D\cap U.
	\end{equation}
	
	Combining these two inequalities with \eqref{nt}, \eqref{appr} and \eqref{imd} provides a constant
	$c_1\in(0,1)$ such that
	$$k_{D\cap V}(z,w)\ge\log\left(1+\frac{g_D(z,w)}{d_z(z)^{1/2}d_z(w)^{1/2}}\right),$$
	where $g_D(z,w)=\max\{c_1|(z-w)_n|-|z-w|^2/c_1,c_1|z-w|^2\}.$ It remains to note that
	$g_D(z,w)\ge\frac{c_1^3}{c_1^2+1}|(z-w)_n|.$
	\smallskip
	
	\noindent{\it Subproof of \eqref{appr}.} Obviously, any translation and unitary transformation
	keeps $\d_D,$ $|u-v|$ and $|u-v|_N.$ So, it suffices to check \eqref{appr} for $\Phi_z=\e_z.$
	
	Observe that $\td_{G_z}\circ\e_z$ is a defining function for $D$ with the same gradient
	as $\d_D$ at $p(z)=0.$ This implies that $d_z=(1+o(1))\d_D$ near $p.$
	
	Further, $z=(z_1,0')$ with $z_1<0$ and then
	$$d_z(w)-d_z(z)=2\Re\langle w-z,\bar\partial d_z(z)\rangle+O(|w-z|^2)$$
	$$=(1+2c_zz_1)\Re(z_1-w_1)+O(|w-z|^2).$$
	
	The proof is finished.
	\smallskip
	
	\noindent{\it Proof of Proposition \ref{c3}.} Having in mind that $D$ is convexifiable near $p,$
	as well as the localizations \eqref{loc}, we may assume that $D$ is strictly convex with
	$\mathcal C^{2,\alpha}$-smooth boundary.

	Let now $c>0$ and $\FF$ be the class of complex geodesics $\vp\in\O(\D,D)$ such that
	$\diam(\vp)\ge c$ and $\d_D(\vp(0))=\max(\d_D\circ\vp)$ (the last one is always achieved
	after a possible reparametrization). By \cite[Proposition 1]{H1} there exists
$c_1>0$ (independent of $c$) such that $\max|\vp'|\le c_1$ and hence
$$|\zz-\ww|\le c_1|\z-\e|\hbox{ for any }\vp\in\FF.$$

On the other hand, by \cite[Proposition 4]{CHL} (see also \eqref{d2}),
	one may find $c_2>0$ such that
	$\d_D(\vp(0))\ge c_2$ for any $\vp\in\FF.$ Then \cite[Corollary 2.1]{Aba}, see also \cite{H2}, provides $c_3>0$ such that
	$$\d_\D\le c_3(\d_D\circ\vp)\hbox{ for any }\vp\in\FF.$$
	
	Since
	\begin{equation}\label{cru}
	k_D(\zz,\ww)=k_\D(\z,\e)\ge\log\left(1+\frac{|\z-\e|}{2\d_\D(\z)^{1/2}\d_\D(\e)^{1/2}}\right),
	\end{equation}
	it follows that
	$$k_D(\zz,\ww)\ge\log\left(1+\frac{|\zz-\ww|}{2c_1c_3\d_D(\zz)^{1/2}\d_D(\ww)^{1/2}}\right).$$
	
	To complete the proof, it remains to apply Theorem \ref{th:gen}.
	\smallskip
	
	\noindent{\it Proof of Proposition \ref{diam}.} The inequality
	$|\vp'(\z)|_N/|\vp'(\z)|\lesssim\diam(\vp)$ follows from Theorem \ref{h2}.
	
	To prove \eqref{d1}, it remains to show that
	$\d_D(\zz)^{1/2}\lesssim\diam(\vp).$ We may assume that $\z=0,$
	$\vp(0)=0$ and $|\vp'(0)|_N/|\vp'(0)|\lesssim \d_D(0)^{1/2}$ (since
$|\vp'(0)|_N/|\vp'(0)|\lesssim\diam(\vp)$, by Theorem \ref{th:gen}).

By \eqref{ma}, we get that
	$|\vp'(0)|\sim\d_D(0)^{1/2}.$ We may assume that $|\vp_1'(0)|\sim\d_D(0)^{1/2}$
($\vp_1$ is the first component of $\vp$). Then the Schwarz lemma implies that
	$$\d_D(0)^{1/2}\sim|\vp_1'(0)|\le\sup_{\D}|\vp_1|<\diam(\vp).$$
	
	Further, the first inequality in \eqref{d2} is a consequence of \eqref{d1}.
	
	Finally, it follows by \cite[Theorem 8]{NO} that there exists $C>0$
	such that the Euclidian length of $\vp(\gamma)$ does not exceed
	$$C\max_{u\in\vp(\gamma)}\left(\d_D(u)^{1/2}\log\frac{1}{\d_D(u)}\right),$$
	for any real geodesic $\gamma$ w.r.t.~$k_\D.$
	This implies second inequality in \eqref{d2} even in the $\mathcal C^2$-smooth case.
	\smallskip
	
	\noindent{\it Proof of Proposition \ref{c33}.} As in Proposition \ref{c3}, we may assume that
	is $D$ is strictly convex with $\mathcal C^{2,\alpha}$-smooth boundary.
	
	Let $\vp\in\O(\D,D)$ be a complex geodesic.
	Since $$|\zz-\ww|\le|\z-\e|\cdot\max|\vp'|$$
	and
	$$\frac{1}{2}\log\frac{\d_D(\vp(0))}{\d_D(\vp(\theta))}\le k_D(\vp(0),\vp(\theta))=k_\D(0,\theta)
	<\frac{1}{2}\log\frac{2}{\d_\D(\theta)},$$ it follows by \eqref{cru} and Theorem \ref{th:gen}  that
	$$k_D(\zz,\ww)\ge\log\left(1+\frac{\d_D(\vp(0))}{4\max|\vp'|}\cdot
	\frac{|\zz-\ww|}{\d_D(\zz)^{1/2}\d_D(\ww)^{1/2}}\right)$$
	$$\ge\log\left(1+\frac{\d_D(\vp(0))}{4C\max|\vp'|.\diam(\vp)}\cdot
	\frac{|(\zz-\ww)_N|}{\d_D(\zz)^{1/2}\d_D(\ww)^{1/2}}\right).$$
	To complete the proof, it remains to apply \eqref{conj}.
	\smallskip
	
\begin{remark} The above proof together with \eqref{d2} shows that $c=c(\eps)$ in Proposition~\ref{c3}
can be chosen such that $c\sim-\eps/\log^2\eps$ for $\eps\in(0,1/2).$
\end{remark}

	\noindent{\it Proof of Corollary \ref{seq}.} (a) follows from Theorem \ref{h1}.
	
	Let us prove (b).  By Theorems \ref{th:gen} and \ref{h2}, we may assume that, after possible reparametrizations,
	$\vp_{z_j;X_j}, \psi_{z_j;w_j}\in\FF$ for some $c>0,$ and $\vp_{z_j;X_j}^{-1}(z_j),
	\psi_{z_j;w_j}^{-1}(z_j)\to 1.$ Note that \cite[Lemma 4]{H1} implies that,
	up to subsequences, $\vp_{z_j;X_j}\to\vp$ and $\psi_{z_j;w_j}\to\psi$ in the $\mathcal C^1$-norm
	on $\ov\D,$ where $\vp$ and $\psi$ are complex geodesics. Clearly,
	$\vp(1)=\psi(1)=p,$ $\vp'(1)=\lambda X$ and $\psi'(1)=\mu X$ for some $\lambda,\mu\in\C.$
	The Hopf lemma implies that $\lambda/\lambda_X,\mu/\lambda_X>0.$ After automorphisms of $\D,$
	we may get $\lambda=\mu=\lambda_X.$

	\section{Proof of Theorem \ref{th:gen}}\label{mt}

	In the proof we can assume that $D$ is strictly convex and $\varphi$ is a geodesic. Actually, the result is trivial when the diameter of an extremal mapping passing through $z$ and $w$ is big. If, in turn, it is small enough, the extremal mapping lies entirely in a neighborhood of some boundary point $p$, which can be exposed globally to a strictly convex one --- see \cite{H1, H2} for details (for another exposing theorem see also \cite[Theorem~1.1]{DFW}). Then this mapping is a stationary disc (see \cite{L2}) and consequently a complex geodesic.

Let us recall here again that the core of the our argument relies on a careful analysis of the so called Lempert theory that was proven in \cite{L1} in $\mathcal C^6$ settings. How to extend it to $\mathcal C^{2,\alpha}$ case was explained in a series of papers by Huang and Lempert, see e.g. discussion \cite[Remark B]{H3} and references contained therein. Below we shall briefly sketch these ideas
\begin{remark}\label{LHremark}

For a strictly linearly convex domain $D$ with $\mathcal C^{2,\alpha}$-smooth boundary, $z\in D$ and $X\in \CC^n$ let $f_{D, z, X}$ denote a complex geodesic in $D$ such that $f_{D,z,X}(0)=z$ and $f'_{D,z,X}(0) = \lambda X$ for some $\lambda>0$.

As noted by Lempert \cite[Lemma 3.1]{L3} (see also \cite{C}), geodesics $f_{D,z,X}$ are $\mathcal C^{1,\alpha}$ smooth on $\TT$. Using the fact that their $\mathcal C^{1/2}$ norms are uniformly bounded when $z$ is within a compact subset of $K$ of $D$ (note that this bound remains uniform when $D$ is perturbed a bit in $\mathcal C^{2}$-topology) Lempert showed in \cite[Proposition 3]{L4} that the bound on $\mathcal C^{1,\alpha}$ norm is in fact locally uniform; see also  \cite[Lemma 4]{H1} for a detailed proof of this fact. Summing up, $\mathcal C^{1,\alpha}$ norms of $f_{D,z,X}$ remain locally uniform with respect to $z$ in $K$ and when $D$ is perturbed in $\mathcal C^{2,\alpha}$ topology.

Therefore one can deduce the following: if $D_n$ converges to a strictly convex domain $D$ in $\mathcal C^{2,\alpha}$ topology, $K$ is a compact subset of $D$ and $f_n$ is a sequence of geodesics in $D_n$ such that $f_n(0)\in K$, then we can find a subsequence converging in $\mathcal C^{1, \alpha}$ topology to a geodesic in $D$. (Using this together with \cite[Theorem 4.1]{L3} and uniqueness of geodesics in $D$ we can even say more: $f_{D,z,X}$ varies continuously, as elements of $\mathcal C^{1,\alpha} (\TT)$, with respect to $z$, $X$ and $D$ varying in $\mathcal C^{2, \alpha}$ topology.)
\end{remark}	
	
We start the proof of Theorem~\ref{th:gen}. Take $\tilde z\in D$ sufficiently close to $\partial D$ and a geodesic $\varphi$ passing through it. The case that is interesting for us is when $\diam(\varphi)\to 0$ (then, of course, $\tilde z\to \partial D$, as geodesics are proper).

We shall use some transformations. Let us mention that they as well as all function that will appear below can be chosen to depend continuously on $\tilde z$, and all constants will be uniformly bounded. Recall the notation from the introduction: if $X\in \CC^n$, then $|X_N|=2|\langle X, g_D(\tilde z)\rangle |$.
	
Translate the domain $D$ so that the closest point to $\tilde z$ on $\partial D$ is $0$. Using unitary map we can make the Levi form at $0$ of a defining function $r$ of $\partial D$ to be a positive diagonal matrix. Replacing the coordinates $z_j$ by $\lambda_j z_j$ we can additionally assume that it is the identity matrix. Using unitary map again we can assume that the tangent space is $\re z_1=0$, i.e. that $r(z) =  2 \re z_1 + |z|^2 + 2 \re  P(z) + O(|z|^{2+\alpha})$, where $P$ is a quadratic polynomial. All these transformations used so far (denote the composition of all of them by $F$) satisfy the following estimate:
	$$
	\frac{|(\tilde z-w)_N|}{|\tilde z-w|} \sim \frac{|(F(\tilde z) - F(w))_N|}{|F(\tilde z)-F(w)|}.
	$$
Moreover, the point $\tilde z$ is mapped a point of the form $(-s,0)$, where $s>0$ (in particular, $(\tilde z - w)_N = \tilde z_1 - w_1$ and $(F(\tilde z) - F(w))_N = F_1(\tilde z) - F_1(w)$). Then composing $F$ with an additional transformation $z_1\mapsto z_1 +P(z)$, $z'\mapsto z'$ gives a biholomorphic map $G$ in a neighborhood $U$of $0$ such that
	  \begin{equation}\label{eq:G1} \frac{|F_1(\tilde z) -F_1(w)|}{|F(\tilde z)-F(w)|} = \frac{|G_1(\tilde z) -G_1(w)|}{|G(\tilde z) - G(w)|} + O(|\tilde z| + |\tilde z-w|)
    \end{equation}	
and the boundary of $G(D\cap U)$ near $0$ is of the form $\{2\re \zeta_1 + |\zeta|^2 + O(|\zeta|^{2+\alpha})<0\}$. To see the estimate in \eqref{eq:G1} it suffices to express $P(\tilde z) - P(w) = \sum (\tilde z_j - w_j) \alpha_j(\tilde z, w)$, where $\alpha_j$ are linear, which implies that $|P(\tilde z) -P(w)|\leq |\tilde z-w| O(|(\tilde z, w)|)$. Thus \eqref{eq:G1} follows from the trivial estimate $|w|\leq |\tilde z| + |\tilde z - w|$.

The translation moves $0$ to $e_1=(1,0')$ and $G(D\cap U)$ to $\{ |\zeta|^2 + O(|\zeta-e_1|^{2+\alpha})<1\}$. Then point $\tilde z$ is mapped to a point of the form $(1-s + \gamma s^2, 0')$, where $\gamma\in \CC$. Clearly $s\to 0$ if $\diam (\varphi)\to 0$. Moreover, $d_D(\tilde z) = s + O(s^2)$.

If $\varphi$ is a geodesic passing through $\tilde z$ and $w$, then $\diam \varphi\gtrsim |\tilde z-e_1|,|\tilde z-w|$ (the first inequality follows from the fact that geodesics are proper and that $|\tilde z-e_1|\sim \d_D(\tilde z)$, the second is trivial). Thus it suffices to focus on achieving estimates for the term $|G_1(\tilde z) - G_1(w)|/|G(\tilde z) - G(w)|$.

\smallskip

Summing up our situation boils down to the following one: $D$ near $e_1$ is of the form $\{|\zeta|^2 + O(|\zeta - e_1|^{2+\alpha})<1\}$, $\tilde z = (1-s + \gamma s^2,0)$, where $s$ is close to $0$ and $\gamma$ is uniformly bounded. We then aim at proving that $|\tilde z_1 - w_1|\leq C|\tilde z-w|\diam( \varphi)$ for any geodesic $\varphi$ in $D$ passing through $\tilde z$ with small diameter and any $w$ in its range.

	Let $m_t(\lambda) = \frac{\lambda + t}{1+ t \lambda},$ $A_t(z) = \left( m_t(z_1), \sqrt{1-t^2}\frac{z'}{1+ t z_1}\right)$, and $r_t(z) = \frac{|1+ tz_1|^2}{1-t^2} r(A_t(z))$, $t\in (0,1)$. Then $r_t$ converges to $\rho(z) =-1+ |z|^2 $ in the $\mathcal C^{2+\alpha}$-topology, as $t\to 1$, when restricted to $\{\re z_1>-1/2\}$  (see \cite[p. 468]{L1}; see also \cite{Kos,KW}).
	
	Take $t$ such that an analytic disc $A_t^{-1}\circ \varphi$ lies entirely in $\{\re z_1>0\}$ and its boundary touches $\{\re z_1=0\}$, say at point $\eta\in \TT$. This, means that
	\begin{equation}\label{eq:t}\frac{t}{1+t^2} = \frac{\re \varphi_1(\eta)}{1 + |\varphi_1(\eta)|^2}.
	\end{equation}
	Note that $t\to 1$ when $\diam(\varphi)\to 0$.
	It is also clear that the above discs are stationary maps, and thus geodesics, in $A_t^{-1} (D) \cap \{\re z_1>-1/2\}$
\smallskip
	
\noindent{\it Claim.} There exists $C_1>0$ such that  $\diam(A_t^{-1}\circ \varphi)>C_1$ for any $t$ as above. In particular, all of
$A_t^{-1}\circ \varphi$ intersect a compact subset $K$ of $\mathbb B\cap \{\re z_1 >-1/2\}$ (compare with the proof of \cite[Proposition 4]{CHL}).
\smallskip
	
	To prove the claim note that it suffices to show that $A_t^{-1}(\varphi(\eta))$ and $A_t^{-1} (\varphi(0))=A_t^{-1}(\tilde z)$ are far away from each other. If this were not the case, $A_t^{-1}(\varphi(\eta))$ would be close to $(\pm i, 0')$, as $A_t^{-1} (\tilde  z)= (m_t^{-1}(1- s+ \gamma s^2), 0').$
	
	Expressing $m_t(a + ib) = 1-s'+i \gamma' s'$, where $a,b,s'$ and $\gamma'$ are real, we see that then $\gamma' = \frac{(1+t)b}{(1+ t a)(1-a) - tb^2}$. In particular, if we write $m_t^{-1}(1- s + \gamma s^2) = a+ ib,$ $a, b\in \RR$, then $b$ must be small (express $1-s+ \gamma s^2 = 1-s' +i O(s'^2)\in \RR+ i \RR$, where $s'= s + O(s^2)$). This proves the claim.
\smallskip
	
We shall apply Remark~\ref{LHremark} to domains $D_t\cap\{\re z_1\geq -1/2\}$ (using cut-off functions we can extend them to globally strictly convex domains). By the claim, after a possible reparametrization, the geodesic $A_t^{-1}\circ \varphi$ is close in $\mathcal C^{1+\alpha-\epsilon}$-topology to a geodesic, say $\psi$, in the ball (recall that $t$ is close to $1$ when $\diam(\varphi)$ is small enough). A geometric characterization of complex geodesics in the Euclidean ball is that they are intersections of the ball with complex affine lines. Since $\psi$ is contained in $\{\re z_1 \geq 0\}$, we get that its range cannot lie on a horizontal line (that is a line of the form $\CC\times \{b\}$). In particular,
$|(\psi(\lambda) - \psi(\mu) )_T|/|\lambda -\mu|$ is a strictly positive constant for $\lambda,\mu\in \overline \DD$.
Here and in the sequel we denote $X_1=(X_1,0)$ and $X_T=(0,X_2, \ldots, X_d)$ . Since $A_t^{-1}\circ \varphi$ is close to $\psi$ in $\mathcal C^1$-topology we easily get that $|(A_t^{-1}\circ \varphi(\lambda) - A_t^{-1}\circ \varphi(\mu) )_T| >D |\lambda -\mu|$, where $D>0$ is a constant uniform with respect to $\lambda, \mu\in \overline{\mathbb D}$. From this we deduce the inequality $$|( A^{-1}_t(\varphi(\lambda))- A_t^{-1}(\varphi(\mu)))_1| <C |(A^{-1}_t(\varphi(\lambda))- A_t^{-1}(\varphi(\mu)))_T|,$$ where $C>0$ is uniform with respect to $\lambda,\mu\in \overline\DD$. Since $\varphi$ passes through $\tilde z$, we thus get that for any $w$ lying on the range of $\varphi$ the following estimate holds:
\begin{equation}\label{eq:est}\frac{|(A^{-1}_t(\tilde z)- A_t^{-1}(w))_1|}{ |( A^{-1}_t(\tilde z)- A_t^{-1}(w))_T|} <C,
\end{equation}
provided that $\diam(\varphi)$ is sufficiently small.

Denoting $x=A^{-1}_t(\tilde z)$, $y= A^{-1}_t(w),$ we rewrite \eqref{eq:est} as
$|x_1 - y_1|/|y_T|<C.$ Let us compute
	\begin{equation*}\label{eq:a}\frac{|(\tilde z-w)_1|}{|(\tilde z-w)_T|} = \frac{|(A_t(x)-A_t(y))_1|}{|(A_t(x) -A_t(t))_T|} = \frac{\sqrt{1-t^2} |x_1-y_1|}{|(1+t x_1)y_T|}\leq C\sqrt{1-t^2}.
	\end{equation*}
	
We shall estimate the right-hand side of the inequality above. It follows from \eqref{eq:t} that $\frac{(1-t)^2}{1+ t^2} = \frac{|1 - \varphi_1(\eta)|^2}{1+ |\varphi_1(\eta)|^2}$, so
$$\sqrt{1-t^2} \sim |1-\varphi_1(\eta)|\leq |\varphi(\eta) - e_1|\leq |\varphi(\eta) - \tilde z| + |\tilde z-e_1|\lesssim \diam (\varphi)$$
(recall that $|\tilde z - e_1|\sim \d_D(\tilde z)$), which completes the proof.

\medskip

\noindent{\bf Acknowledgements.} The authors are grateful to Pascal J. Thomas for careful reading of the manuscript. We would also like to thank the anonymous referee for numerous remarks that substantially improved the paper.

\end{document}